\documentclass[a4paper, 12pt]{article}
\usepackage{amssymb, amsmath, amsfonts}
\usepackage{amsthm,amscd,amsfonts,latexsym,color,epsfig}
\begin{document}

\begin{center}
{\bf A BOUNDARY PROBLEM WITH INTEGRAL GLUING CONDITION FOR A PARABOLIC-HYPERBOLIC EQUATION INVOLVING THE CAPUTO FRACTIONAL DERIVATIVE}\\
\medskip
Erkinjon T. Karimov , Jasurjon S. Akhatov\\

\end{center}

\bigskip

\textbf{Abstract. }In the present work we investigate the Tricomi problem with an integral gluing condition for parabolic-hyperbolic equation involving the Caputo fractional differential operator. Using the method of energy integrals we prove the uniqueness of the solution for considered problem. The existence of the solution have been proved applying methods of ordinary differential equations and Fredholm integral equations. Solution is represented in an explicit form.

\bigskip

\bigskip

\section{Formulation of the problem}

It is well-known that investigations of fractional analogs of main ODE and PDEs one motivated by their appearance in building mathematical models for real-life processes [1-3]. They are as well interesting for mathematicians as natural generalizations of integer order ODE and PDEs. Specialists in the theory of boundary problems for PDEs began to develop it in this direction. There are many works [4-6] devoted to the investigation of various boundary problems for PDEs.

A distinctive side of this work is the usage of gluing condition of the integral form, containing regular continuous gluing condition as a particular case. We note that for the first time boundary problem with integral gluing condition for parabolic-hyperbolic type equation was used in the work [7]. Then some generalizations of this work were published in [8-9]. Special gluing condition of the integral form for parabolic-hyperbolic equation with the Riemann-Liouville fractional differential operator was discussed in [10].

In the present work we use integral gluing condition with kernel, which has a more general form than the kernel used in [11]. When we prove the uniqueness of the solution we must put some restrictions to the kernel (see theorem 1), but for the existence of solution we don't need these conditions (see theorem 2).

Consider the equation
$$
0=\left\{ \begin{array}{l}
u_{xx}-{}_{C}D_{0y}^{\alpha }u,\hfill y>0, \\
u_{xx}-u_{yy},\hfill y<0 \\
\end{array}
\right.
\eqno (1)
$$
in the domain $\Omega =\Omega^+\cup\Omega^-\cup AB$,  where $0<\alpha <1$, $AB=\left\{ \left( x,y \right):\,0<x<1,\,y=0 \right\},$
$\Omega^+=\left\{ \left( x,y \right):\,0<x<1,\,0<y<1 \right\},$ $\Omega^-=\left\{ \left( x,y \right):\,-y<x<y+1,\,-1/2<y<0 \right\},$\\
$$
{}_{C}D_{0y}^{\alpha }f=\frac{1}{\Gamma \left( 1-\alpha  \right)}\int\limits_{0}^{y}{{{\left( y-t \right)}^{-\alpha }}{f}'\left( t \right)dt}
$$
is the Caputo fractional differential operator of order $\alpha \,\left( 0<\alpha <1 \right)$, $\Gamma \left( \cdot  \right)$ is the Euler's gamma-function [12].

\textbf{Problem.} Find a solution of the equation (1) belonging to
$$
W=\left\{ u\left( x,y \right):\,u\in C\left( \overline{\Omega } \right)\cap {{C}^{2}}\left( {{\Omega }^{-}} \right),\,{{u}_{xx}}\in C\left( {{\Omega }^{+}} \right),\,{}_{C}D_{0y}^{\alpha }u\in C\left( {{\Omega }^{+}} \right)\, \right\}
$$
satisfying the boundary conditions
$$
u\left( 0,y \right)={{\varphi }_{1}}\left( y \right),\,\,\,0\le y\le 1,\eqno (2)
$$
$$
u\left( 1,y \right)={{\varphi }_{2}}\left( y \right),\,\,\,0\le y\le 1,\eqno (3)
$$
$$
u\left( x,-x \right)=\psi \left( x \right),\,\,\,0\le x\le 1/2, \eqno (4)
$$
and the gluing condition
$$
\underset{y\to +0}{\mathop{\lim }}\,{{y}^{1-\alpha }}{{u}_{y}}\left( x,y \right)={{\gamma }_{1}}{{u}_{y}}\left( x,-0 \right)+{{\gamma }_{2}}\int\limits_{0}^{x}{{{u}_{y}}\left( t,-0 \right)Q\left( x,t \right)dt},\,\,0<x<1.\eqno (5)
$$
Here ${{\varphi }_{1}},\,{{\varphi }_{2}},\,\psi ,\,Q\left( \cdot ,\cdot  \right)$ are given functions, such that
${{\varphi }_{1}}\left( 0 \right)=\psi \left( 0 \right),\,\,\,{{\gamma }_{1}},{{\gamma }_{2}}=const,\,\gamma _{1}^{2}+\gamma _{2}^{2}\ne 0$.

\section{Uniqueness of the solution}

Let us set
$$
u\left( x,+0 \right)={{\tau }_{1}}\left( x \right),\,0\le x\le 1,\,\,u\left( x,-0 \right)={{\tau }_{2}}\left( x \right),\,0\le x\le 1,
$$
$$
{{u}_{y}}\left( x,-0 \right)={{\nu }_{2}}\left( x \right),\,\,0<x<1,\,\,\underset{y\to +0}{\mathop{\lim }}\,{{y}^{1-\alpha }}{{u}_{y}}\left( x,y \right)={{\nu }_{1}}\left( x \right),\,0<x<1,
$$
$$
{{u}_{x}}\left( x,+0 \right)={{{\tau }'}_{1}}\left( x \right),\,0<x<1,\,\,{{u}_{x}}\left( x,-0 \right)={{{\tau }'}_{2}}\left( x \right),\,0<x<1.
$$

It is known that solution of the Cauchy problem for Eq.(1) in $\Omega^-$ can be represented as follows:
$$
u\left( x,y \right)=\frac{1}{2}\left[ {{\tau }_{1}}\left( x+y \right)+{{\tau }_{2}}\left( x-y \right)-\int\limits_{x-y}^{x+y}{{{\nu }_{2}}\left( t \right)dt} \right].\eqno (6)
$$
After using condition (4) in (6) we find
$$
\tau'_2\left( x \right)-2\psi'\left( x/2 \right)=\nu_2\left( x \right),\,\,0<x<1.\eqno (7)
$$
From the Eq.(1) at $y\rightarrow +0$ we get [13]
$$
\tau''_1\left( x \right)-\Gamma \left( \alpha  \right)\nu_1\left( x \right)=0. \eqno (8)
$$

Below we prove the uniqueness of the solution of the formulated problem. For this aim, first we suppose that the problem has two solutions, then denoting the difference of these solutions by $u$ we will get appropriate homogeneous problem. If we prove that this homogeneous problem has only the trivial solution, then we can state that the original problem has unique solution.

We multiply equation (8) by $\tau_1(x)$ and integrate from $0$ to $1$:
$$
\int\limits_0^1{\tau''_1\left( x \right)\tau_1\left( x \right)dx}-\Gamma \left( \alpha  \right)\int\limits_0^1{\tau_1\left( x \right)\nu_1\left( x \right)dx}=0.\eqno (9)
$$
We investigate the integral $I=\int\limits_0^1{\tau_1\left( x \right)\nu_1\left( x \right)dx}$. Considering the gluing condition (5), we have
$$
\nu_1\left( x \right)=\gamma_1\nu_2\left( x \right)+\gamma_2\int\limits_0^x{\nu_2\left( t \right)Q\left( x,t \right)dt},\,\,0<x<1.\eqno (10)
$$
In the homogeneous case, i.e. $\psi \left( x \right)=0$, from (7) we obtain $\nu_2\left( x \right)=\tau'_2\left( x \right)$, hence (10) will be written as
$$
\nu_1\left( x \right)=\gamma_1\tau'_2\left( x \right)+\gamma_2\int\limits_0^x{\tau'_2\left( t \right)Q\left( x,t \right)dt},\,\,0<x<1.\eqno (11)
$$
We substitute (11)  into the integral $I$ and consider $\tau_1\left( 0 \right)=0,\,\tau_1\left( 1 \right)=0$ (which are deduced from conditions (2), (3) in the homogeneous case), we have
$$
I=\int\limits_0^1{\tau_1\left( x \right)\nu_1\left( x \right)dx}=\gamma_2\int\limits_0^1{\tau_1^2\left( x \right)Q\left( x,x \right)dx}-\gamma _2\int\limits_0^1{\tau_1\left( x \right)dx\int\limits_0^x{\tau_1\left( t \right)\frac{\partial }{\partial t}Q\left( x,t \right)dt}}.\eqno (12)
$$
Let
$$
\frac{\partial }{\partial t}Q\left( x,t \right)=-Q_1\left( x \right)Q_1\left( t \right).
$$
Then
$$
I=\gamma_2\int\limits_0^1{\tau_1^2\left( x \right)Q\left( x,x \right)dx}+\frac{\gamma_2\Phi^2\left( 1 \right)}{2},\eqno (13)
$$
where
$$
\Phi \left( x \right)=\int\limits_0^x{\tau_1\left( t \right)Q_1\left( t \right)dt},\,
Q\left( x,t \right)=Q\left( x,0 \right)-\int\limits_0^t{Q_1\left( x \right)Q_1\left( z \right)dz}.
$$
From (9) and (13) we get
$$
\int\limits_0^1{{\tau'}_1^2\left( x \right)dx}+\Gamma \left( \alpha  \right)\gamma_2\left[ \int\limits_0^1{\tau_1^2\left( x \right)Q\left( x,x \right)dx}+\frac{\Phi^2\left( 1 \right)}{2} \right]=0. \eqno (14)
$$
Since $\Gamma(\alpha)>0$ for $0<\alpha<1$, then if $\gamma_2\ge 0,\,\,Q\left( x,x \right)>0$ from (14) we easily get $\tau_1(x)=0$ for any $x\in \left[ 0,1 \right]$.

Based on the solution of the first boundary problem for Eq.(1) in the domain $\Omega^+$ we obtain $u(x,y)\equiv 0$ in $\overline{\Omega^+}$.
Since $u(x,y) \in C\left(\overline{\Omega}\right)$, we get that $u(x,y)\equiv 0$ in $\overline{\Omega}$.

\smallskip

Hence, we proved the following

\smallskip

\textbf{Theorem 1.} Let $\gamma_2\ge 0,\,\frac{\partial }{\partial t}Q\left( x,t \right)=-Q_1\left( x \right)Q_1\left( t \right)$ and $Q\left( x,x \right)>0$. If there exists a solution to problem, then it is unique.

Example of $Q(x,t)$:
$$
Q\left( x,t \right)=e^{-x}\left( 1+e^{-t} \right).
$$

\section{Existence of the solution}
From (7), (8) and (10) we have
$$
{\tau''}_1(x)-A\tau_1(x)=F_1(x), \eqno (15)
$$
where $A=\Gamma(\alpha)\gamma_1,$
$$
F_1(x)=\gamma_2\Gamma(\alpha)\int\limits_0^x{{\tau'}_1(t)Q(x,t)dt}-\Gamma(\alpha)\left[\gamma_1\psi\left(\frac{x}{2}\right)
+\gamma_2\int\limits_0^x{\psi'\left(\frac{t}{2}\right)Q(x,t)dt}\right].\eqno (16)
$$
The solution of the equation (15) together with the conditions
$$
\tau_1(0)=\psi(0),\,\,\tau_1(1)=\varphi_2(0) \eqno (17)
$$
has the form
$$
\tau_1(x)=\frac{1}{1-e^A}\left[\varphi_2(0)\left(1-e^{Ax}\right)+\psi(0)\left(e^{Ax}-e^A\right)\right]+\int\limits_0^1{G_0(x,\xi)F_1(\xi)d\xi},\eqno (18)
$$
where
$$
G_0(x,\xi)=\frac{1}{A\left[e^{Ax}-e^{A(x-1)}\right]}
\left\{
\begin{array}{l}
\left(1-e^{A\xi}\right)\left(1-e^{A(x-1)}\right),\,0\leq \xi \leq x,\\
\left(1-e^{A(\xi-1)}\right)\left(1-e^{Ax}\right),\,x\leq \xi \leq 1\\
\end{array}
\right.\eqno (19)
$$
is the Green's function of the problem (15), (17).
Considering (16) and integrating by parts, we obtain
$$
\tau_1(x)-\int\limits_0^1{\tau_1(\xi)K(x,\xi)d\xi=F_2(x)}, \eqno (20)
$$
where
$$
K(x,\xi)=\gamma_2\Gamma(\alpha)\left[G_0(x,\xi)Q(\xi,\xi)+\int\limits_{\xi}^1{G_0(\xi,t)\frac{\partial}{\partial \xi}Q(t,\xi)dt}\right], \eqno (21)
$$
$$
F_2(x)=\frac{1}{1-e^A}\left[\varphi_2(0)\left(1-e^{Ax}\right)+\psi(0)\left(e^{Ax}-e^A\right)\right]-
$$
$$
-\Gamma(\alpha)\int\limits_0^1{G_0(x,\xi)\left[\gamma_1\psi\left(\frac{\xi}{2}\right)
+\gamma_2\int\limits_0^{\xi}{\psi'\left(\frac{t}{2}\right)Q(\xi,t)dt}\right]d\xi}.\eqno (22)
$$
Since the kernel $K(x,\xi)$ is continuous and $F_2(x)$ is continuously differentiable, the solution of integral equation (20) can be written via the resolvent-kernel:
$$
\tau_1(x)=F_2(x)-\int\limits_0^1{F_2(\xi)R(x,\xi)d\xi},\eqno (23)
$$
where $R(x,\xi)$ is the resolvent-kernel of $K(x,\xi)$.

\smallskip

The unknown functions $\nu_1(x)$ and $\nu_2(x)$ will be found by the following formulas:
$$
\nu_1(x)=\frac{1}{\Gamma(\alpha)}\left[F_2''(x)-\int\limits_0^1{F_2(\xi)\frac{\partial^2}{\partial x^2}R(x,\xi)d\xi}\right],
$$

$$
\nu_2(x)=F_2'(x)-\int\limits_0^1{F_2(\xi)\frac{\partial}{\partial x}R(x,\xi)d\xi}-\psi'\left(\frac{x}{2}\right).
$$
The solution of the problem in the domain $\Omega^+$ can be written as follows
$$
u(x,y)=\int\limits_0^y{G_{\xi}(x,y,0,\eta)\varphi_1(\eta)d\eta}-\int\limits_0^y{G_{\xi}(x,y,1,\eta)\varphi_2(\eta)d\eta}+
\int\limits_0^1{\overline{G}(x-\xi,y)\tau_1(\xi)d\xi},\eqno (24)
$$
where
$$
\overline{G}(x-\xi,y)=\frac{1}{\Gamma(1-\alpha)}\int\limits_0^y{\eta^{-\alpha}G(x,y,\xi,\eta)d\eta},
$$
$$
G(x,y,\xi,\eta)=\frac{(y-\eta)^{\beta-1}}{2}\sum\limits_{n=-\infty}^{\infty}\left[e_{1,\beta}^{1,\beta}\left(-\frac{|x-\xi+2n|}{(y-\eta)^\beta}\right)-
e_{1,\beta}^{1,\beta}\left(-\frac{|x+\xi+2n|}{(y-\eta)^\beta}\right)\right]
$$
is the Green's function of the first boundary problem for Eq.(1) in the domain $\Omega^+$ with the Riemann-Liouville fractional differential operator instead of the Caputo ones [13], $\beta=\alpha/2$,
$$
e_{1,\beta}^{1,\beta}(z)=\sum\limits_{n=0}^{\infty}\frac{z^n}{n!\Gamma(\beta-\beta n)}
$$
is the Wright type function [13].

The solution of the problem in the domain $\Omega^-$ will be found by the formula (6).
Hence, we proved the following

\textbf{Theorem 2.} If $\varphi_i(y),\psi(x)\in C[0,1]\cap C^1(0,1)$, $Q(x,t)\in C^1\left([0,1]\times[0,1]\right)$, then there exists a solution of the problem and it can be represented in the domain $\Omega^+$ by formula (24) and in the domain $\Omega^-$ by the formula (6).

\section{Acknowledgement}
Authors grateful for Professor M.Kirane for his useful suggestions, which made paper more readable.

\smallskip

\centerline{\bf References}
\begin{enumerate}
\item {\it C.G.Koh and J.M.Kelly.} Application of fractional derivatives to seismic analysis of base-isolated models. Earthquake Engineering and Structural Dynamics, 19 (1990), pp.229-241.\\
\item{\it T.J.Anastasio.} The fractional order dynamics of brainstem vestibulo-oculomotor neurons. Biological Cybernetics, 72 (1994), pp. 69-79.\\
\item{\it A.A.Kilbas, O.A.Repin.} An analog of the Tricomi problem for a mixed type equation with a partial fractional derivative. Fractional Calculus and Applied Analysis, 13(1) (2010), pp.69-84.\\
\item{\it Y.Luchko.} Initial-boundary-value problems for the generalized multi-term time-fractional diffusion equation. Journal of Mathematical Analysis and Applications, 374 (2011), pp. 538-548.\\
\item{\it A.V.Pskhu.} Solution of the first boundary value problem for a fractional order diffusion equation. Differential Equations, 39(9) (2003), pp.1359-1363.\\
\item{\it Y.Povstenko.} Neumann boundary-value problems for a time-fractional diffusion-wave equation in a half-plane. Computer and Mathematics with Applications, doi:10.1016/j.camwa.2012.02.064\\
\item{\it N.Yu.Kapustin and E.I.Moiseev.} On spectral problems with a spectral parameter in the boundary condition. Differential Equations, 33(1) (1997), pp.116-120.\\
\item{\it B.E.Eshmatov and E.T.Karimov.} Boundary value problems with continuous and special gluing conditions for parabolic-hyperbolic type equations. Central European Journal of Mathematics, 5(4) (2007), pp.741-750.
\item{\it E.T.Karimov.} Non-local problems with special gluing condition for the parabolic-hyperbolic type equation with complex spectral parameter. PanAmerican Mathematical Journal, 17 (2007), pp.11-20.
\item{\it A.S.Berdyshev, A.Cabada and E.T.Karimov.} On a non-local boundary problem for a parabolic-hyperbolic equation involving a Riemann-Liouville fractional differential operator. Nonlinear Analysis, 75 (2012), pp.3268-3273.\\
\item{\it A.S.Berdyshev, E.T.Karimov and N.Akhtaeva.} Boundary value problems with integral gluing conditions for fractional-order mixed-type equation. International Journal of differential Equations, (2011), Article ID 268465.\\
\item {\it I. Podlubny.} Fractional Differential Equations, Academic Press, San Diego, 1999.\\
\item {\it A.V.Pskhu.} Solution of boundary value problems for the fractional diffusion equation by the Green function method. Differential Equations, 39(10) (2003), pp.1509-1513.\\

\end{enumerate}

Erkinjon T. Karimov\\
\medskip
E-mail: erkinjon@gmail.com\\
\emph{Institute of Mathematics, National University of Uzbekistan named after Mirzo Ulugbek (Tashkent, Uzbekistan),}\\

\bigskip

Jasurjon S. Akhatov\\
\medskip
E-mail: jahatov@gmail.com\\
\emph{Physical-Technical Institute, SPA "Physics-Sun", Academy of Sciences of the Republic of Uzbekistan (Tashkent, Uzbekistan)}\\

\end{document}